\documentclass[12pt]{article}

\usepackage{amsmath}
\usepackage{amssymb}
\usepackage{amsthm}
\usepackage{bbm}
\usepackage[dvips]{graphicx}
\usepackage{color}
\usepackage{epsfig}
\usepackage[mathscr]{eucal}

\newtheorem{theorem}{Theorem}

\newtheorem{prop}[theorem]{Proposition}
\newtheorem{lemma}[theorem]{Lemma}
\newtheorem{cor}[theorem]{Corollary}

\newtheorem{question}{Question}
\newtheorem{conjecture}[question]{Conjecture}
\newtheorem{problem}[question]{Problem}
\newtheorem{conj}[question]{Conjecture}

\newcommand{\ft}{\ensuremath{\mathbb{F}_q}}

\newcommand{\R}{{\mathbb R}}

\newcommand{\by}{{\bf y}}
\newcommand{\bx}{{\bf x}}

\newcommand{\cC}{{\mathcal C}}

\newcommand{\cL}{{\mathcal L}}
\newcommand{\cM}{{\mathcal M}}

\allowdisplaybreaks

\title{Collinear Triple Hypergraphs and the Finite Plane Kakeya Problem}
\author{Joshua Cooper} 

\begin{document}

\maketitle

\begin{abstract} We show that the problem of counting collinear points in a permutation (previously considered by the author and J. Solymosi in \cite{CS05}) and the well-known finite plane Kakeya problem are intimately connected.  Via counting arguments and by studying the hypergraph of collinear triples we show a new lower bound ($5q/14 + O(1)$)  for the number of collinear triples of a permutation of $\ft$ and a new lower bound ($q(q+1)/2 + 5q/14 + O(1)$) on the size of the smallest Besicovitch set in $\ft^2$.  Several interesting questions about the structure of the collinear triple hypergraph are presented.
\end{abstract} 

\section{Introduction}

The Kakeya problem, solved in 1928 by Besicovitch, asks for the smallest subset of $\R^2$ so that a unit-length interval may be continuously rotated $360^\circ$ within it.  The striking answer is that there are sets of arbitrarily small measure that permit such a ``full turn.''  If one does not require that the interval be continuously rotated, but simply that the set contain a unit interval in every direction, there are even sets of measure zero that do the job.  This statement also holds in higher dimensions.  Perhaps unsurprisingly, such sets must still be ``large,'' i.e., have Hausdorff (or Minkowski) dimension equal to that of the ambient space.  To show this is the case in dimensions higher than two is a very important open problem, with deep connections to additive number theory and harmonic analysis.

In his influential manuscript \cite{W99}, Wolff describes the ``finite field Kakeya problem'' and makes a convincing argument that the problems that arise are very similar to those of the analytic setting.  The finite field Kakeya problem can be stated as follows.  A {\it Besicovitch set} in $\ft^d$ is a subset that contains a line in every direction, i.e., a translate of every line through the origin.  What is the smallest Besicovitch set?  The natural conjecture is that it has cardinality $c q^d$, but this is not known in dimensions greater than 2.  (Indeed, the best result for dimension 3 is only $q^{5/2 + 10^{-10}}$ \cite{KLT00}.)

It is easy to see that, for the finite plane, every Besicovitch set has cardinality at least $q^2/2$.  In fact, this result appears already in \cite{W99}.  On the other hand, there is a natural construction that gives an upper bound of $q(q+1)/2 + (q-1)/2$ for $q$ odd.  Hence, we have the following conjecture.

\begin{conjecture} For $q$ an odd prime power, the smallest Besicovitch set in $\ft^2$ has cardinality $q(q+1)/2 + (q-1)/2$.
\end{conjecture} 
(The situation when $q$ is a power of two is considerably simpler.)  Recently, Faber \cite{F06} provided a very elegant argument that, in fact, the smallest finite plane Besicovitch set has cardinality at least $q(q+1)/2 + q/3$, getting quite close to the upper bound.

In the current paper, we show that (1) Faber's argument dualizes to improve the best known bound for the number of collinear triples of a permutation (see \cite{CS05} and the next section); and (2) with some additional work, it is possible to prove the following theorem:

\begin{theorem} \label{thm:mainthm} The smallest Besicovitch set in $\ft^2$ has cardinality at least
$$
\frac{q(q+1)}{2} + \frac{5q}{14} - \frac{1}{14}.
$$
\end{theorem}

Note that $5/14 = 1/3 + 1/42$, so this is indeed an improvement.  In fact, this lower bound gives the conjectural upper bound for $q = 3,5,7,9$.  We include below several very interesting problems that arise from our analysis.

\section{Preliminaries}

Let $q$ be an odd prime power.  Write $\ft$ for the finite field of order $q$ and $\ft^\ast$ for its multiplicative group, i.e., $\ft \setminus \{0\}$.  Let $f : \ft \rightarrow \ft$ be some function; then we define the ``graph of $f$'', denoted $\Gamma_f$, to be the set $\{(j,f(j)) : j \in \ft\} \subset \ft^2$.  We will write $m_{ij} = m_{ij}(f)$ for the slope of the line connecting the points in $\Gamma_f$ corresponding to the ordinates $i$ and $j$, i.e., $m_{ij} = (f(j) - f(i))/(j-i)$.  We refer to the slope of a collinear set of points in the obvious sense, and call two collinear sets with the same slope ``parallel.''  In general, a slope is an element of $PG(1,q)$, and therefore can take on the value $\infty$ if the denominator of this expression is $0$, although this is clearly never the case for a graph $\Gamma_f$.  Also, the slope of points in $\Gamma_\sigma$ cannot be $0$ for $\sigma$ a permutation.   Finally, we write $\ell(\bx,m)$ for the line in $\ft^2$ through $\bx$ with slope $m$.

Suppose that $P$ is a minimal Besicovitch set in $\ft^2$.  Clearly, we may assume that $P$ is of the form
\begin{equation} \label{eq:minimalBesicovitch}
P = \bigcup_{m \in PG(1,q)} \ell(\bx_m,m).
\end{equation}
Write $\mu(\by) = |\{m \in PG(1,q) : \by \in \ell(\bx_m,m) \}|$, i.e., the multiplicity of $\by$ in the representation of $P$ by (\ref{eq:minimalBesicovitch}).  We have the following ``incidence formula'' for the size of $P$ from \cite{F06}:

\begin{theorem}[Incidence Formula, \cite{F06}] \label{thm:incidenceformula}
$$
|P| = \frac{q(q+1)}{2} + \sum_{\by \in P} \binom{\mu(\by)-1}{2}.
$$
\end{theorem}

Note that $\binom{\mu(\by)}{3} \geq \binom{\mu(\by)-1}{2} \geq \mu(\by) - 2 \geq 1$ for $\mu(\by) \geq 3$.  The following result (actually, an equivalent version) appears in \cite{F06}.  For a Besicovitch set $P$ of the form (\ref{eq:minimalBesicovitch}), define $\rho(m)$ to be the sum of $\mu(\by)-2$ over all $\by \in \ell(\bx_m,m)$ with $\mu(\by) \geq 3$.

\begin{prop} \label{prop:onlyonezero} Given a Besicovitch set $P = \bigcup_{m \in PG(1,q)} \ell(\bx_m,m)$, at most one slope $m \in PG(1,q)$ has $\rho(m) = 0$.
\end{prop}

It is not hard to conclude from this Proposition and Theorem \ref{thm:incidenceformula} that $|P| \geq q(q+1)/2 + q/3$.  We can go further with this result, however.  Indeed, suppose that $\min_{m \in PG(1,q)} \rho(m) = R$.  If $R \geq 2$, then 
\begin{align*}
|P| & = \frac{q(q+1)}{2} + \sum_{\by \in P} \binom{\mu(\by)-1}{2} \\
& = \frac{q(q+1)}{2} + \sum_{m \in PG(1,q)} \sum_{\substack{\by \in \ell(\bx_m,m) \\ \mu(\by) \geq 3}} \frac{1}{\mu(\by)} \binom{\mu(\by)-1}{2} \\
& = \frac{q(q+1)}{2} + \sum_{m \in PG(1,q)} \sum_{\substack{\by \in \ell(\bx_m,m) \\ \mu(\by) \geq 3}} \left ( \frac{1}{2} - \frac{1}{2 \mu(\by)} \right ) \cdot (\mu(\by)-2) \\
& \geq \frac{q(q+1)}{2} + \sum_{m \in PG(1,q)} \sum_{\substack{\by \in \ell(\bx_m,m) \\ \mu(\by) \geq 3}} \frac{\mu(\by) - 2}{3} \\
& = \frac{q(q+1)}{2} + \sum_{m \in PG(1,q)} \frac{\rho(m)}{3} \\
& \geq \frac{q(q+1)}{2} + \frac{2(q+1)}{3}.
\end{align*}
Hence, Theorem \ref{thm:mainthm} follows unless $R = 0$ or $R = 1$.

Now, we may apply a nonsingular affine transformation of $\ft^2$ to make the line $\ell(\bx_m,m)$ with $\rho(m) = R$ become the vertical line through the origin while retaining the property of $P$ being a minimal Besicovitch set.  Dualizing the {\it nonvertical} lines of $P$ under the map $\{(x,y):y = mx+b\} \mapsto (m,b)$ gives a set $S \subset \ft^2$ which is the graph $\Gamma_\sigma$ of some function $\sigma$.  If $R = 0$, then $\sigma$ is a permutation; if $R = 1$, then $\sigma$ is a ``semipermutation'', that is, a function whose range has cardinality $q-1$.  If $R = 0$, points $\by$ with $\mu(\by) \geq 3$ correspond to collinear $\mu(\by)$-tuples in $\Gamma_\sigma$; if $R = 1$, then all $\by$ with $\mu(\by) \geq 3$ except the single point of multiplicity $\geq 3$ on $\ell((0,0),\infty)$  correspond to collinear $\mu(\by)$-tuples in $\Gamma_\sigma$.

This leads us to the following definition.  Given $S \subset \ft^2 = \ft \times \ft$, define the 3-uniform hypergraph $\cC(S)$ on the vertex set $S$ to be the set of all collinear triples $\{\bx_1 , \bx_2 , \bx_3\} \subset S$.  Let $\cL(S)$ for the hypergraph on the vertex set $S$ whose edges are the maximal collinear subsets of $S$.  For a hypergraph $H$, write $\|H\|$ for the quantity $\sum_{e \in E(H)} (|e|-1)(|e|-2)/2$; we write $\|S\| = \|\cL(S)\|$ for convenience.  Note that the number of edges in $\cC(S)$, which we write simply $|\cC(S)|$, is at least $\|S\|$.  It is clear by the Incidence Formula that, when $R = 0$ and $\sigma$ is the permutation arising from dualizing $P$,
$$
|P| = \frac{q(q+1)}{2} + \|\Gamma_\sigma\|.
$$
When $R = 1$ and $\sigma$ is the semipermutation arising from dualizing $P$,
$$
|P| = 1 + \frac{q(q+1)}{2} + \|\Gamma_\sigma\|.
$$
Our problem therefore reduces to that of determining the minimum value of $\|\Gamma_\sigma\|$ when $\sigma$ is a permutation or a semipermutation.  In both cases, we may assume that every point of $\Gamma_\sigma$ appears in at least one edge of $\cC(\Gamma_\sigma)$ by the minimality of $R$ and Proposition \ref{prop:onlyonezero}.  In fact, we reprove this statement below for permutations, since the argument is so short and elegant.
 
The following conjecture appears in \cite{CS05}, along with a proof of the lower bound $(q-1)/4$.

\begin{conj} For any permutation $\sigma: \ft \rightarrow \ft$, the number of collinear triples in $\Gamma_\sigma$  is at least $(q-1)/2.$
\end{conj}

If this conjecture is true, it is best possible, as evidenced by the graph of the function $s \mapsto s^{-1}$ for $s \neq 0$ and $0 \mapsto 0$.  Unfortunately, we still cannot show the full conjecture.  However, we do show below that every permutation $\sigma$ satisfies $\|\Gamma_\sigma\| \geq (5q-1)/14$, from which it follows that the number of collinear triples in $\Gamma_\sigma$ is at least this quantity.  Furthermore, we show that every semipermutation $\sigma$ so that $\cC(\Gamma_\sigma)$ has no isolated points satisfies $\|\Gamma_\sigma\| \geq 5(q-1)/14$.  In the former case, then, $|P| \geq q(q+1)/2 + (5q-1)/14$, and in the latter,
$$
|P| \geq \frac{q(q+1)}{2} + \frac{5q-5}{14} + 1 = \frac{q(q+1)}{2} + \frac{5q}{14} + \frac{9}{14}.
$$
Theorem \ref{thm:mainthm} follows.

\section{Permutations}

Throughout this section we assume that $\sigma : \ft \rightarrow \ft$ is a permutation.  The next proposition states that no point of $\Gamma_\sigma$ is isolated in $\cC(\Gamma_\sigma)$, and the argument is essentially the dual of that appearing in \cite{F06}.

\begin{prop} \label{prop:nolonely} For every $j \in \ft$, there exists an $E \in \cC(\Gamma_\sigma)$ so that $(j,\sigma(j)) \in E$.
\end{prop}
\begin{proof} Clearly we may assume that $(j,\sigma(j))=(0,0)$.  Suppose $(j,\sigma(j))$ were an isolated point in $\Gamma_\sigma$, i.e., there exists no $E \in \cC(\Gamma_\sigma)$ so that $(j,\sigma(j)) \in E$.  Consider the slopes $m_{0i}$ for $i \neq 0$, and, in particular, the product
$$
M = \prod_{i \in \ft^\ast} m_{0i}.
$$
Since $(0,0)$ does not participate in any collinear triple, the line $\ell((0,0),m)$ contains at most two points for each $m \in \ft^\ast$ (including the origin).  However, since there are $q-1$ such lines and $q-1$ points $(j,\sigma(j))$ for $j \neq 0$, each such line contains exactly two points of $\Gamma_\sigma$, by the pigeonhole principle.  This means that $M$ is the product of all elements of $\ft^\ast$, which it is easy to see is $-1$.  On the other hand,
$$
M = \prod_{i \in \ft^\ast} \frac{\sigma(i)}{i} = \frac{ \prod_{i \in \ft^\ast} i}{ \prod_{i \in \ft^\ast} i} = 1,
$$
a contradiction.
\end{proof}

In fact, Proposition \ref{prop:nolonely} can be pushed considerably further.

\begin{prop} Suppose that $\bx=(j,\sigma(j))$ participates in no collinear quadruples.  Denote by $\cM_i$, $i = 1,2,3$, the set of slopes $m \in \ft^\ast$ so that $\ell(\bx;m)$ contains exactly $i$ points of $\Gamma_\sigma$.  Then $|\cM_1| = |\cM_3|$, and
$$
\prod_{m \in \cM_1} m = -\prod_{m \in \cM_3} m.
$$
\end{prop}
\begin{proof} First of all, $|\cM_2| + 2 |\cM_3| = q-1$, since there are $q-1$ points in $\Gamma_\sigma \setminus \{\bx\}$, and $|\cM_1| + |\cM_2| + |\cM_3| = q-1$, since there are $q-1$ slopes of lines through $(j,\sigma(j))$ which are not vertical or horizontal.  Subtracting the two equations, we see that $|\cM_1| = |\cM_3|$.

We may assume without loss of generality that $\bx = (0,0)$.  Consider again the product
$$
M = \prod_{j \in \ft^\ast} \frac{\sigma(j)}{j}.
$$
On the one hand, the numerator and denominator of this expression both enumerate all elements of $\ft^\ast$, so that $M = 1$.  On the other hand,
\begin{align*}
M & = \left ( \prod_{m \in \cM_2} m \right ) \left ( \prod_{m \in \cM_3} m^2 \right ) \\
& = \left ( \prod_{m \in \ft^\ast} m \right ) \left ( \prod_{m \in \cM_1} m \right )^{-1} \left ( \prod_{m \in \cM_3} m \right ) \\
& = - \left ( \prod_{m \in \cM_1} m \right )^{-1} \left ( \prod_{m \in \cM_3} m \right ).
\end{align*}
Hence,
$$
\prod_{m \in \cM_1} m = -\prod_{m \in \cM_3} m.
$$
\end{proof}

Call a point $p \in \Gamma_\sigma$ which belongs to exactly one collinear triple ``lonely,'' i.e., 
$$
|\{E \in \cC(\Gamma_\sigma) : p \in E\}| = 1.
$$

\begin{prop} \label{prop:perpline} If $\bx = (j,\sigma)$ is lonely, then
\begin{enumerate}
\item There is one slope $m_j \in \ft^\ast$ so that $\ell(\bx,m_j)$ contains three points of $\Gamma_\sigma$, one slope $m^\perp_j$ so that $\ell(\bx,m^\perp_j)$ contains one point of $\Gamma_\sigma$, and all other slopes $m \neq m_j,m^\perp_j$ are such that $\ell(\bx,m)$ contains two points.
\item $m_j^\perp = -m_j$.
\end{enumerate}
\end{prop}
\begin{proof} Since $\bx$ participates in exactly one collinear triple, it is easy to see that, letting $m_j$ be the slope of that triple, $|\ell(\bx,m_j) \cap \Gamma_\sigma| = 3$ and $|\ell(\bx,m) \cap \Gamma_\sigma| \leq 2$ for all other $m \in \ft^\ast \setminus \{m_j\}$.  There are $q-2$ other slopes, while there are only $q-3$ other points $(i,\sigma(i))$.  Hence, only $q-3$ of the lines $\ell(\bx,m)$, $m \in \ft^\ast \setminus \{m_j\}$ actually have a point of $\Gamma_\sigma$ besides $\bx$ on them, leaving one line with no other point.  This slope we call $m_j^\perp$.

To see the second claim, we need only apply the previous proposition with $\cM_3 = \{m_j\}$, $\cM_1 = \{m_j^\perp\}$.
\end{proof}

\begin{problem} Can an edge of $\cC(\Gamma_\sigma)$ be isolated, i.e., its three points lonely?
\end{problem}

If not, then we may conclude that $\|\Gamma_\sigma\| \geq 2q/5$ by a simple counting argument.  Indeed, we have the following.

\begin{prop} \label{prop:twofifths} If $T \subset \Gamma_\sigma$ of cardinality $t$ induces no isolated edges and no isolated points in $\cC(\Gamma_\sigma)$, then
$$
\left \| \Gamma_\sigma \cap T \right \| \geq \frac{2 t}{5}.
$$
\end{prop}
\begin{proof} Clearly, it suffices to prove a lower bound of $2/5$ on the quotient $\|H\|/t$ for all connected hypergraphs $H$ on $t$ vertices with at least $2$ edges and so that each edge contains at least three points.  Note that, for any such $H$, it is possible to order its edges $e_1,\ldots,e_k$ so that each $e_j$, $2 \leq j \leq k$, intersects $e_i$ for some $i < j$.  Define
$$
a_j = \left | e_j \setminus \bigcup_{i < j} e_i \right |.
$$
Then $a_1 = |e_1|$ and $a_j \leq |e_j| - 1$ for $j > 1$.  Since $t = \sum_{j =1}^k a_j$,
\begin{align*}
\frac{\|H\|}{t} &= \frac{\sum_{e \in H} (|e|-1)(|e|-2)}{2 \sum_{j =1}^k a_j} \\
& \geq \frac{\sum_{e \in H} |e|^2 - 3 \sum_{e \in H} |e| + 2 |H|}{2 (\sum_{e \in H} |e| - |H| + 1)}\\
& \geq \frac{|H|^{-1}(\sum_{e \in H} |e|)^2 - 3 \sum_{e \in H} |e| + 2 |H|}{2 (\sum_{e \in H} |e| - |H| + 1)}
\end{align*}
by Cauchy-Schwarz.  Let $\epsilon = \sum_{e \in H} |e|$ and $h = |H|$.  Then this reduces to
$$
\frac{\|H\|}{t} \geq \frac{\epsilon^2 - 3 \epsilon h+ 2 h^2}{2 h (\epsilon - h + 1)}
$$
which it is easily checked over $h \geq 2$, $\epsilon \geq 3h$ attains its minimum of $2/5$ at $h = 2$, $e = 6$.
\end{proof}

Define an {\it isolated matching} of a hypergraph $H$ to be a set of (disjoint) isolated edges, i.e., a matching which is not intersected by any other edges of $H$.

\begin{lemma} \label{lemma:oneptlines} For every isolated matching $M \subset \cC(\Gamma_\sigma)$, there are at least $2 |M| + t - \|\Gamma_\sigma\| - 1$ slopes $m \in \ft^\ast$ so that there exists an $\ell(\bx,m)$ whose intersection with $\Gamma_\sigma$ consists of exactly one point of $T$, where $t = |T|$.
\end{lemma}
\begin{proof} Suppose that $M \subset \cC(\Gamma_\sigma)$ has exactly $(q-t)/3$ edges, leaving a set $T = \Gamma_\sigma \setminus \bigcup M$ of $t$ other vertices.  Let $S_1$ be the set of slopes of all collinear triples contained in $T$, so that $|S_1| \leq \|T\|$.  Note that every edge of $\cC(\Gamma_\sigma)$ containing a point of $T$ in fact contains three points of $T$, since $M$ is an isolated matching.  Consider the set of slopes $S_2$ which occur in $M$, or whose negative is an element of $M$.  Then $|S_2| \leq 2 |M| = 2(q-t)/3$, and $|S_1 \cup S_2| \leq \|T\| + 2(q-t)/3$.  Hence, there exist at least
$$
|\ft^\ast \setminus (S_1 \cup S_2)| \geq q - 1 - \left ( \|T\| + \frac{2(q-t)}{3} \right ) = \frac{q+2t}{3} - \|T\| - 1,
$$
slopes $m \in \ft^\ast$ which are not elements of $S_1 \cup S_2$.  Since $ \|T\| = \|\Gamma_\sigma\| - (q-t)/3$, this quantity is at least
\begin{align*}
\frac{q+2t}{3} - \left ( \|\Gamma_\sigma \| - \frac{q-t}{3} \right ) - 1 &= \frac{2q + t}{3} - \|\Gamma_\sigma\| - 1 \\
&= 2 |M| + t - \|\Gamma_\sigma\| - 1.
\end{align*}

If a line $\ell$ of slope $m$ intersects $M$ in a point $\bx$, then, since $m \not \in S_2$ and $M$ is an isolated matching, $\ell$ cannot contain more than one other point of $\Gamma_\sigma$.  Furthermore, since $-m \not \in S_2$, $\ell$ contains at least one other point of $\Gamma_\sigma$.  Now, if a line $\ell$ of slope $m$ does not intersect $M$ but {\it does} intersect $T$, then it cannot intersect $T$ in three or more points, since $m \not \in S_1$.  Altogether, this means that every line of slope $m$ intersects $\Gamma_\sigma$ either (a) nowhere, (b) at two points of $M$, (c) at one point of $T$, (d) at two points of $T$, or (e) at one point of $M$ and one point of $T$.

The conclusion follows from the observation that the quantities $|\ell(\bx,m) \cap \Gamma_\sigma|$ of types (a),(b),(d), and (e) are even, while $q$, which is odd, satisfies
$$
q = \sum_{\ell \textrm{ of slope } m} |\ell \cap \Gamma_\sigma|.
$$
Hence there is an odd number of lines of type (c).
\end{proof}

Now, for a point $\bx \in \Gamma_\sigma$, define $\|\bx\|$ to be the sum over all maximal collinear subsets $\{\by_1,\ldots,\by_k\} \subset \Gamma_\sigma$ of the quantity $\phi(k) = (k-1)(k-2)/(2k)$.

\begin{lemma} \label{lemma:emptymeansfull} Suppose that there are $B$ distinct lines $\ell(\bx,m)$, $m \in \ft^\ast$, through a point $\bx \in \Gamma_\sigma$ so that
$$
\ell(\bx,m) \cap \Gamma_\sigma = \bx.
$$
Then $\| \bx \| \geq B/3$.
\end{lemma}
\begin{proof} Suppose there are $A$ lines $\ell(\bx,m)$ so that $|\ell(\bx,m) \cap \Gamma_\sigma| \geq 3$, and let $a_1,\ldots,a_A$ denote the quantities $|\ell(\bx,m) \cap \Gamma_\sigma|$.  Among the slopes $m \in \ft^\ast$, there are $B$ lines $\ell(\bx,m)$ so that $|\ell(\bx,m) \cap \Gamma_\sigma| = 1$, and $q - 1 - (A+B)$ lines so that $|\ell(\bx,m) \cap \Gamma_\sigma| = 2$.  On the other hand, since $| \Gamma_\sigma | = q$,
$$
q-1-(A+B) + 2A \leq q-1-(A+B) + \sum_{j=1}^A (a_j - 1) = |\Gamma_\sigma \setminus \{\bx\}| =  q - 1,
$$
whence $A \leq B$ and $\sum_{j=1}^A a_j = 2A + B$.  Let $\delta = B/A \geq 1$.  Then, by the convexity and monotonicity of $\phi$ on $[3,\infty)$,
\begin{align*}
\| \bx \| & = \sum_{i=1}^A \phi(a_i) \geq A \phi(\sum_{i=1}^A a_i) = A \phi \left ( \frac{2A + B}{A} \right ) = \frac{B}{\delta} \phi ( 2 + \delta ) \\
& = \frac{B}{\delta} \frac{\delta(\delta + 1)}{2 (\delta + 2)} = \frac{B(\delta + 1)}{2 (\delta + 2)} \geq \frac{B}{3}.
\end{align*}
\end{proof}

\begin{cor} \label{cor:onethird} Suppose that $M \subset \cC(\Gamma_\sigma)$ is an isolated matching, $T = \Gamma_\sigma \setminus \bigcup M$, and there are $R$ distinct lines $\ell$ intersecting $T$ so that $|\ell \cap \Gamma_\sigma| = 1$.  Then 
$$
\| T \| \geq R/3.
$$
\end{cor}
\begin{proof} It is easy to see that, since $M$ is an isolated matching,
$$
\| T \| = \sum_{\bx \in T} \| \bx \|.
$$
The conclusion then follows immediately from Lemma \ref{lemma:emptymeansfull}.
\end{proof}

\begin{theorem} $\| \Gamma_\sigma \| \geq \frac{5q-1}{14}$.
\end{theorem}
\begin{proof} Suppose that $M \subset \cC(\Gamma_\sigma)$ is a maximal isolated matching.  Define $T = \Gamma_\sigma \setminus \bigcup M$ and $t = |T|$.  By Lemma \ref{lemma:oneptlines}, there are at least $2 |M| + t - \| \Gamma_\sigma \| - 1$ distinct lines which intersect $T$ in at least one point.  Therefore, by Corollary \ref{cor:onethird}, $\| T \| \geq  2 |M|/3 + t/3 - \| \Gamma_\sigma \|/3 - 1/3$.  Then
$$
\| \Gamma_\sigma \| = \|T\| + |M| \geq \frac{5 |M|}{3} + \frac{t}{3} - \frac{\| \Gamma_\sigma \|}{3} - \frac{1}{3},
$$
whence
$$
\frac{q-t}{3} = |M| \leq \frac{4 \| \Gamma_\sigma \|}{5} - \frac{t}{5} + \frac{1}{5}
$$
from which we may deduce
$$
t \geq \frac{5q-3}{2} - 6\| \Gamma_\sigma \|.
$$
Since $M$ is maximal, $T$ does not induce any isolated edges.  Therefore, by Proposition \ref{prop:twofifths},
$$
\| \Gamma_\sigma \| = |M| + \|T\| \geq \frac{q-t}{3} + \frac{2t}{5} = \frac{q}{3} + \frac{t}{15} \geq \frac{q}{2} - \frac{2 \| \Gamma_\sigma \|}{5} - \frac{1}{10}.
$$
Rearranging, we find that
$$
\| \Gamma_\sigma \| \geq \frac{5q-1}{14}.
$$
\end{proof}

\section{Semipermutations}

Throughout this section, we assume that $\sigma$ is a function from $\ft$ to itself which repeats exactly one value, i.e., $|\sigma(\ft)| = q-1$, and so that there are no isolated points in $\cC(\Gamma_\sigma)$.  Fix the following notation: $\sigma(\ft) = \ft \setminus \{b\}$, $\sigma(z_1)=\sigma(z_2)=a$ for some $z_1,z_2 \in \ft$, $z_1 \neq z_2$.

\begin{prop} If $\bx = (j,\sigma)$ is lonely for $j \neq z_1,z_2$, then
\begin{enumerate}
\item There is one slope $m_j \in \ft^\ast$ so that $\ell(\bx,m_j)$ contains three points of $\Gamma_\sigma$, one slope $m^\perp_j$ so that $\ell(\bx,m^\perp_j)$ contains one point of $\Gamma_\sigma$, and all other slopes $m \neq m_j,m^\perp_j$ are such that $\ell(\bx,m)$ contains two points.
\item $m_j^\perp = -m_j \cdot (\sigma(j) - b)/(\sigma(j) - a)$.
\end{enumerate}
\end{prop}
\begin{proof} The proof of (1) is identical to that of Proposition \ref{prop:perpline}.  To see the second claim, note that
$$
\prod_{i \neq j} \frac{\sigma(j)-\sigma(i)}{j-i} = \prod_{m \in \ft^\ast} m \cdot \frac{m_j}{m_j^\perp} = -\frac{m_j}{m_j^\perp}.
$$
On the other hand, the numerator of the left-hand side is given by
$$
\prod_{i \neq j} (\sigma(j)-\sigma(i)) = \prod_{s \in \ft^\ast}{s} \cdot \frac{\sigma(j) - a}{\sigma(j) - b} = - \frac{\sigma(j) - a}{\sigma(j) - b},
$$
and the denominator is $\prod_{s \in \ft^\ast} s = -1$, whence
$$
-\frac{m_j}{m_j^\perp} = \frac{\sigma(j) - a}{\sigma(j) - b},
$$
and the conclusion follows.
\end{proof}

\begin{prop} If $\bx = (j,a)$ is lonely for $j = z_1$ (or $z_2$), then
\begin{enumerate}
\item There is one slope $m_j \in \ft^\ast$ so that $\ell(\bx,m_j)$ contains three points of $\Gamma_\sigma$, two slopes $m^\perp_{j1}, m^\perp_{j2}$ so that each $\ell(\bx,m^\perp_{ji})$ contains one point of $\Gamma_\sigma$, and all other slopes $m \neq m,m^\perp_{j1},m^\perp_{j2}$ are such that $\ell(\bx,m_j)$ contains two points.
\item $m^\perp_{j1} m^\perp_{j2} = m_j \cdot (a - b)/(z_2 - z_1)$ (respectively, $m_j \cdot (a - b)/(z_1 - z_2)$).
\end{enumerate}
\end{prop}
\begin{proof} Write $j^\prime = z_2$ (respectively, $z_1$).  Of the $q-1$ lines $\ell(\bx,m)$ through $\bx$, $m \in \ft^\ast$, there is exactly one $m_j$ for which $|\ell(\bx,m_j) \cap \Gamma_\sigma| = 3$ (since $\bx$ is lonely), and all other intersections contain $1$ or $2$ points.  Since $q-2$ of the points of $\Gamma_\sigma$ are counted exactly once by these quantities -- that is, all of $\Gamma_\sigma$ except $\bx$ and $(j^\prime,a)$ -- we have
$$
q-2 = \sum_{m \in \ft^\ast} \left ( |\ell(\bx,m_j) \cap \Gamma_\sigma| - 1 \right ) = 2 + |\{m : |\ell(\bx,m_j) \cap \Gamma_\sigma| = 2 \}|.
$$
This, in turn, implies that there are $q-4$ slopes $m \in \ft^\ast$ so that $\ell(\bx,m_j) \cap \Gamma_\sigma$ consists of exactly two points, and the remaining $(q-1) - (q-4) - 1 = 2$ directions $m^\perp_{j1}, m^\perp_{j2}$ contain only one point: $\bx$.

For the second claim, note that
$$
\prod_{i \neq j,j^\prime} \frac{\sigma(j)-\sigma(i)}{j-i} = \prod_{m \in \ft^\ast} m \cdot \frac{m_j}{m^\perp_{j1} m^\perp_{j2}} = - \frac{m_j}{m^\perp_{j1} m^\perp_{j2}}.
$$
On the other hand, the numerator of the left-hand side is given by
$$
\prod_{i \neq j,j^\prime} (\sigma(j)-\sigma(i)) = \prod_{s \in \ft^\ast}{s} \cdot \frac{1}{a - b} = \frac{1}{b-a},
$$
and the denominator is $\prod_{s \in \ft^\ast} s / (j - j^\prime) = 1/(j^\prime - j)$, whence
$$
\frac{m_j}{m^\perp_{j1} m^\perp_{j2}} = \frac{j^\prime - j}{a - b},
$$
and the conclusion follows.
\end{proof}

\begin{lemma} \label{lemma:oneptlinesagain} For every isolated matching $M \subset \cC(\Gamma_\sigma)$, there are at least $2 |M| + t - \|\Gamma_\sigma\| - 3$ slopes $m \in \ft^\ast$ so that there exists an $\ell(\bx,m)$ whose intersection with $\Gamma_\sigma$ consists of exactly one point of $T$, where $t=|T|$.
\end{lemma}
\begin{proof} Suppose that $M \subset \cC(\Gamma_\sigma)$ has exactly $(q-t)/3$ edges, leaving a set $T = \Gamma_\sigma \setminus \bigcup M$ of $t$ other vertices.  Let $S_1$ be the set of slopes of collinear triples in $T$, so that $|S_1| \leq \|T\|$.  Note that every edge of $\cC(\Gamma_\sigma)$ containing a point of $T$ in fact contains three points of $T$, since $M$ is an isolated matching.  Consider the set $S_2$ of slopes $m_j$ and $m_j^\perp$ which occur in $M$.  Then $|S_2| \leq 2 |M| = 2(q-t)/3 + 2$ (the $2$ arising in case $z_1$ or $z_2 \in \bigcup M$), and $|S_1 \cup S_2| \leq \|T\| + 2(q-t)/3 + 2$.  Hence, there exist at least
$$
|\ft^\ast \setminus (S_1 \cup S_2)| \geq q - 1 - \left (\|T\| + \frac{2(q-t)}{3} + 2 \right ) = \frac{q+2t}{3} - \|T\| - 3,
$$
slopes $m \in \ft^\ast$ which are not elements of $S_1 \cup S_2$.  Since $\|T\| = \|\Gamma_\sigma\| - (q-t)/3$, this quantity is at least
\begin{align*}
\frac{q+2t}{3} - \left (\|\Gamma_\sigma\| - \frac{q-t}{3} \right ) - 3 &= \frac{2q + t}{3} - \|\Gamma_\sigma\| - 3 \\
&= 2 |M| + t - \|\Gamma_\sigma\| - 3.
\end{align*}

If a line $\ell$ of slope $m$ intersects $M$ in a point $\bx$, then, since $m \not \in S_2$ and $M$ is an isolated matching, $\ell$ cannot contain more than one other point of $\cC(\Gamma_\sigma)$.  Furthermore, since $-m \not \in S_2$, $\ell$ contains at least one other point of $\Gamma_\sigma$.  Now, if a line $\ell$ of slope $m$ does not intersect $M$ but {\it does} intersect $T$, then it cannot intersect $T$ in three or more points, since $m \not \in S_1$.  Altogether, this means that every line of slope $m$ intersects $\Gamma_\sigma$ either (a) nowhere, (b) at two points of $M$, (c) at one point of $T$, (d) at two points of $T$, or (e) at one point of $M$ and one point of $T$.

The conclusion follows from the observation that the quantities $|\ell(\bx,m) \cap \Gamma_\sigma|$ of types (a),(b),(d), and (e) are even, while $q$, which is odd, satisfies
$$
q = \sum_{\ell \textrm{ of slope } m} |\ell \cap \Gamma_\sigma|.
$$
Hence there is an odd number of lines of type (c).
\end{proof}

\begin{lemma} \label{lemma:emptymeansfullagain} Suppose that there are $B$ distinct lines $\ell(\bx,m)$, $m \in \ft^\ast$, through a point $\bx \in \Gamma_\sigma$ so that
$$
\ell(\bx,m) \cap \Gamma_\sigma = \bx.
$$
If $\bx \neq z_1,z_2$, then $\| \bx \| \geq B/3$; otherwise $\| \bx \| \geq (B-1)/3$.
\end{lemma}
\begin{proof} If $\bx \neq z_1,z_2$, the proof is identical to that of Lemma \ref{lemma:emptymeansfull}.  So, suppose that $\bx = z_1$ or $z_2$, and that there are $A$ lines $\ell(\bx,m)$ so that $|\ell(\bx,m) \cap \Gamma_\sigma| \geq 3$.  Let $a_1,\ldots,a_A$ be the quantities $|\ell(\bx,m) \cap \Gamma_\sigma|$.  Among the slopes $m \in \ft^\ast$, there are $B$ lines $\ell(\bx,m)$ so that $|\ell(\bx,m) \cap \Gamma_\sigma| = 1$, and $q - 1 - (A+B)$ lines so that $|\ell(\bx,m) \cap \Gamma_\sigma| = 2$.  On the other hand, since $| \Gamma_\sigma | = q$,
$$
q-1-(A+B) + 2A \leq q-1 - (A+B) + \sum_{j=1}^A (a_j - 1) = |\Gamma_\sigma \setminus \{z_1,z_2\}| =  q - 2,
$$
whence $A \leq B - 1$ and $\sum_{j=1}^A (a_j - 1) = 2A+B-1$.  Let $\delta = (B-1)/A \geq 1$.  Again, by the convexity and monotonicity of $\phi$ on $[3,\infty)$,
\begin{align*}
\| \bx \| & = \sum_{i=1}^A \phi(a_i) \geq A \phi(\sum_{i=1}^A a_i) = A \phi \left ( \frac{2A + B - 1}{A} \right ) = \frac{B - 1}{\delta} \phi \left ( 2 + \delta \right) \\
& = \frac{B-1}{\delta} \cdot \frac{\delta(\delta + 1)}{2 (\delta + 2)} = \frac{(B-1)(\delta + 1)}{2 (\delta + 2)} \geq \frac{B-1}{3}.
\end{align*}
\end{proof}

\begin{cor} \label{cor:onethirdagain} Suppose that $M \subset \cC(\Gamma_\sigma)$ is an isolated matching, $T = \Gamma_\sigma \setminus \bigcup M$, and there are $R$ distinct lines $\ell$ intersecting $T$ so that $|\ell \cap \Gamma_\sigma| = 1$.  Then $\| T\| \geq (R-2)/3$.
\end{cor}
\begin{proof} Since $M$ is an isolated matching,
$$
\| T \| = \sum_{\bx \in T} \| \bx \|.
$$
The conclusion then follows immediately from Lemma \ref{lemma:emptymeansfullagain}, since at most two of the $\bx \in T$ are $z_1$ or $z_2$.
\end{proof}

\begin{theorem} $\| \Gamma_\sigma \| \geq \frac{5q-5}{14}$.
\end{theorem}
\begin{proof} Suppose that $M \subset \cC(\Gamma_\sigma)$ is a maximal isolated matching.  Define $T = \Gamma_\sigma \setminus \bigcup M$ and $t = |T|$.  By Lemma \ref{lemma:oneptlinesagain}, there are at least $2 |M| + t - \|\Gamma_\sigma\| - 3$ distinct lines which intersect $T$ in at least one point.  Then, by Corollary \ref{cor:onethirdagain}, $\| T \| \geq 2 |M|/3 + t/3 - \|\Gamma_\sigma\|/3 - 5/3$, from which it follows that 
$$
\| \Gamma_\sigma \| = |M| + \| T \| \geq \frac{5 |M|}{3} + \frac{t}{3} - \frac{\| \Gamma_\sigma \|}{3} - \frac{5}{3},
$$
whence
$$
\frac{q-t}{3} = |M| \leq \frac{4 \| \Gamma_\sigma \|}{5} - \frac{t}{5} + 1
$$
from which we may deduce
$$
t \geq \frac{5q-15}{2} - 6\| \Gamma_\sigma \|.
$$
Since $M$ is maximal, $T$ does not induce any isolated edges.  Therefore, by Proposition \ref{prop:twofifths},
$$
\| \Gamma_\sigma \| = |M| + \| T \| \geq \frac{q-t}{3} + \frac{2t}{5} = \frac{q}{3} + \frac{t}{15} \geq \frac{5q}{12} - \frac{\| \Gamma_\sigma \|}{5} - \frac{1}{4}.
$$
Rearranging, we find that
$$
\| \Gamma_\sigma \| \geq \frac{5q-5}{14}.
$$
\end{proof}

\section{Questions}

First of all, although it is clear that an understanding of the size of the maximum isolated matching in $\cC(\Gamma_\sigma)$ would improve our results, we have made little progress in determining what its extremal values are.  The arguments in previous sections can be exploited to achieve the following relatively weak bound on $t$.

\begin{prop} Let $\sigma : \ft \rightarrow \ft$ be a permutation.  The largest isolated matching $M$ in $\cC(\Gamma_\sigma)$ has size less than $(q-\sqrt{q/7})/3$.
\end{prop}
\begin{proof} The number of collinear triples in $\Gamma_\sigma$ is, as we have seen, at least $(5q-1)/14$.  If there were $t$ points in $T = \Gamma_\sigma \setminus \bigcup M$, the number of collinear triples spanned by these points is at most $\binom{t}{2}/3 = t(t-1)/6$.  Hence,
$$
\frac{q-t}{3} + \frac{t(t-1)}{6} \geq \frac{5q-1}{14}.
$$
Solving for $t$ yields
$$
t \geq \frac{3}{2} + \sqrt{\frac{q}{7} + \frac{51}{28}} > \sqrt{\frac{q}{7}}.
$$
\end{proof}

As mentioned before, we do not even know if it is possible to have a single isolated edge in $\cC(\Gamma_\sigma)$.  Of course, the biggest outstanding problem at this time is closing the remaining $q/7$ gap both in the minimum number of collinear triples of a permutation and the size of the smallest Besicovitch set in $\ft^2$.

\section{Acknowledgements}

The author wishes to thank Fan Chung of UCSD and Angelika Steger of ETH for their generous support during the period this paper was produced.  Thanks also to Sebi Cioab\v{a} and Kevin Costello for valuable discussions, and to Xander Faber for his indispensable help in whipping this paper into shape.

\end{document}